\input amstex
\documentstyle{amsppt}
\topmatter
\magnification=\magstep1
\pagewidth{5.2 in}
\pageheight{6.7 in}
%\hcorrection{-0.4in}
%\vcorrection{-0.4in}
\abovedisplayskip=12pt \belowdisplayskip=12pt
\NoBlackBoxes
\title A note on $q$-Volkenborn Integration \endtitle

%Use \endgraf to indicate new paragraph
%\thanks will become a 1st page footnote
\author  Taekyun Kim \endauthor
\affil{{\it Institute of Science Education}\\
{\it Kongju National University, Kongju 314-701, S. Korea }\\
{\rm e-mail : tkim$\@$kongju.ac.kr}}\\\\
\endaffil

\define \bb{{\beta}}

\define\C{\Bbb C_p}

\define\Z{\Bbb Z_p}

\define\pv{\par \vskip0.2cm}

\abstract{In this paper,  we construct $q$-numbers and polynomials
by using the $q$-Volkenborn integrals.
 } \endabstract
 \subjclass  11S80 \endsubjclass
\thanks  \endthanks
\leftheadtext{T. Kim}
\rightheadtext{$q$-Volkenborn Integration }
%\TagsOnRight
\endtopmatter

\document

\head 1. Introduction \endhead Let $p$ be a fixed odd prime, and
let $\Bbb C_p$ denote the $p$-adic completion of the algebraic
closure of $\Bbb Q_p.$ For $d$ a fixed positive integer with
$(p,d)=1$, let
$$X=X_d=\varprojlim_N \Bbb Z/dp^N , \;\;X_1=\Bbb Z_p,$$
$$X^*=\bigcup\Sb 0<a<dp\\ (a,p)=1\endSb a+dp\Bbb Z_p,$$
$$a+dp^N\Bbb Z_p=\{x\in X\mid x\equiv a\pmod{dp^N}\},$$
where $a\in \Bbb Z$ lies in $0\leq a<dp^N ,$ (cf. [1], [2]).

The $p$-adic absolute value in $\Bbb C_p$ is normalized so that
$|p|_p=\frac1p .$ Let $q$ be variously considered as an
indeterminate a complex number $q \in \Bbb C$, or a $p$-adic
number $q\in\Bbb C_p .$ If $q\in\Bbb C ,$ we always assume
$|q|<1.$ If $q\in\Bbb C_p,$ we always assume $|q-1|_p <
p^{-\frac{1}{p-1}}$, so that $q^x= \exp(x\log q)$ for $|x|_p \le
1.$ Throughout this paper, we use the following notation :
$$[x]_q=[x:q]= \frac{1-q^x}{1-q}.$$
We say that $f$ is uniformly differentiable function at a point
 $ a \in \Bbb Z_p$-- and denote this property by $ f \in UD(\Bbb
 Z_p )$-- if the difference quotients
 $$ F_f (x, y)=\frac{f(x)-f(y)}{x-y} ,$$ have a limit $l=
 f^{\prime}(a)$ as $( x, y) \rightarrow (a,a) ,$ cf. [1].
For $f\in U D(\Bbb Z_p),$ let us start with the expression
$$\frac1{[p^N]_q} \sum_{0 \le j < p^N} q^i f(j) =\sum_{0\le j<p^N} f(j)
\mu_q(j+p^N\Bbb Z_p), \text{ cf. [2], }$$ representing
$q$-analogue of Riemann sums for $f$.

The integral of $f$ on $\Bbb Z_p$ will be defined as limit
($n\rightarrow \infty$) of these sums, when it exists. The
$q$-Volkenborn integral of a function $f \in {\text{UD}} (\Bbb
Z_p) $ is defined by
$$\int_{\Bbb Z_p} f(x) d\mu_q (x)= \lim_{N\rightarrow \infty} \frac1{[p^N]_q}
\sum_{0\le j<p^N}
f(j)q^j.$$
Note that if $f_n \rightarrow f$ in $UD( Z_p)$; then
$$\int_{\Bbb Z_p} f_n(x) d \mu_q(x) \rightarrow \int_{\Bbb Z_p} f(x) d
\mu_q(x).$$ In [1, 2], the $q$-Bernoulli polynomials are defined
by
$$ \beta_{m, q}(x)=
\int_{\Bbb Z_p}[ x+x_1]_q^m d\mu_q(x_1).$$
 In this paper, we investigate some properties of  $q$-Volkenborn integrals
 on $\Bbb Z_p$ in the meaning of fermionic and give
some formulae which are related to $q$-numbers and polynomials in
these integrals.

 \head \S 2. On the applications of $q$-Volkenborn integrals \endhead

For any positive integer $N,$ we set
$$\mu_q(a+dp^N\Z)=\frac{q^a}{[dp^N]_q},$$
and this can be extended to distribution on $X$ (see [1]).\pv This
distribution yields an integral for each non-negative integer $m$:
$$\align
\int_{\Z}[a]_q^m\,d\mu_q(a)&=\int_X[a]_q^m\,d\mu_q(a)\\
&=\bb_{m,q}=\frac{1}{(1-q)^m}\sum_{i=0}^m\pmatrix m\\i\endpmatrix
(-1)^i\frac{i+1}{[i+1]}.
\endalign$$
\indent The $q$-Bernoulli polynomials in variable $x$ in $\C$ with
$|x|_p\leq 1$ are defined by
$$\bb_{n,q}(x)=\int_{\Z}[x+t]_q^nd\mu_q(t).$$
These can be written as
$$\align
\bb_{n,q}(x)&=\sum_{i=0}^n\pmatrix n\\i\endpmatrix q^{ix}\bb_{i,q}[x]_q^{n-i}\\
&=\frac{1}{(1-q)^m}\sum_{i=0}^n\pmatrix n\\i\endpmatrix (-1)^i
q^{xi}\frac{i+1}{[i+1]_q}, \ (\text{cf.}\ [3]).\endalign $$ In the
meaning of fermionic, we consider the $q$-numbers by using
$q$-Volkenborn integrals as follows: $$\int_{X_m}
[x]_q^kd\mu_{-q}(x)=\int_{\Bbb Z_p}[x]_q^k d\mu_{-q}(x)=K_{k,q}
\text{ for $k, m\in\Bbb N$ }.$$ Thus, we note that
$$K_{k,q}=[2]_q\left(\frac{1}{1-q}\right)^k\sum_{l=0}^k\binom kl
(-1)^l\frac{1}{1+q^{l+1}}, \tag1$$ where $\binom ki$ is the
binomial coefficient.

For $x\in\Bbb Z_p,$ we define $q$-polynomials as follows:
$$\int_{\Bbb Z_p}[x+y]_q^k d\mu_{-q}(y)=K_{k,q}(x). \tag2$$
By (2), we easily see that
$$K_{k,q}(x)=\sum_{n=0}^k\binom kn [x]_q^{k-n}q^{nx}K_{n,q}. \tag3$$
We now observe that
$$\frac{1}{[p^N]_{-q}}\sum_{n=0}^{p^N-1}[x+y]_q^n(-q)^y
=[2]_q\left(\frac{1}{1-q}\right)^n\sum_{k=0}^n\binom nk
(-1)^kq^{xk}\frac{1}{1+q^{p^N}}\frac{1+q^{p^N(k+1)}}{1+q^{k+1}}.\tag4$$
Thus, we have
$$K_{n,q}(x)=\int_{\Bbb Z_p}[x+y]_{q}^n
d\mu_{-q}(y)=[2]_q\left(\frac{1}{1-q}\right)^n\sum_{k=0}^n\binom
nk(-1)^kq^{xk}\frac{1}{1+q^{k+1}}.\tag5$$ By using the definition
of Eq.(2), we will give the distribution of $q$-polynomials. From
the definition of $q$-Volkenborn integral, we derive the below
formula:
$$\int_{X_m }[x+y]_q^n d\mu_{-q}(y)
=\frac{[m]_q^m}{[m]_{-q}}\sum_{a=0}^{m-1}(-1)^aq^a\int_{\Bbb
Z_p}[\frac{a+x}{m}+y]_{q^m}^n d\mu_{-q^m}(y), \text{ if $m$ is
odd}.\tag6$$ By (2) and (6), we obtain the below equation. Let $m$
be the odd integer. Then we have
$$K_{n,q}(x)=\frac{[m]_q^n}{[m]_{-q}}\sum_{a=0}^{m-1}(-1)^aq^aK_{n,q^m}(\frac{a+x}{m}).\tag7$$
Let $\chi$ be the primitive Dirichlet's character with conductor
$f\in\Bbb N.$ Then we also consider the below generalized
$q$-Bernoulli numbers attached to $\chi$:
$$K_{n,\chi,q}=\int_{X_f}\chi(x)[x]_q^n d\mu_{-q}(x). \tag8$$
By simple calculation in Eq.(8), we easily see that
$$K_{n,\chi,q}=\frac{[f]_q^n}{[f]_{-q}}\sum_{a=0}^{f-1}\chi(a)(-1)^aq^a
\int_{\Bbb Z_p}[\frac{a}{f}+x]_{q^f}^nd\mu_{-q^f}(x), \text{ if
$f$ is odd integer.}\tag 9$$ By (9) and (2), we obtain the
following:
$$K_{n,\chi,q}=\frac{[f]_q^n}{[f]_{-q}}\sum_{a=0}^{f-1}\chi(a)(-1)^aq^aK_{n,q^f}(\frac{a}{f}),\tag10$$
where $\chi$ is a primitive Dirichlet's character with conductor
odd integer $f$. It seems to be interesting to study the Kummer
congruences for the generalized $q$-numbers attached to $\chi$ in
the Eq.(10) because the $q$-Volkenborn measure $d\mu_{-q}(x)\sim
\frac{[2]_q}{2}(-1)^xq^x $.

\head \S 3. Conclusion \endhead
 In complex plane, we assume that
$q\in\Bbb C$ with $|q|<1 .$ It was well known that the Euler
numbers are defined by
$$\frac{2}{e^t+1}=\sum_{n=0}^{\infty} E_n \frac{t^n}{n!}, \text{
 for $|t|<2\pi$, cf. [1, 2, 3].}\tag11$$
 Let $F_q(t)=\sum_{n=0}^{\infty}K_{n,q}\frac{t^n}{n!}$. Then we
 have
 $$F_{q}(t)=e^{\frac{t}{1-q}}\sum_{j=0}^{\infty}\frac{1+q}{1+q^{j+1}}(-1)^j\big(\frac{1}{1-q}\big)^j\frac{t^j}{j!}.
 \tag12$$
By (12), we can find the nice generating function of $K_{n,q}$ as
follows:
$$F_{q}(t)=\sum_{n=0}^{\infty}K_{n,q}\frac{t^n}{n!}=[2]_q\sum_{n=0}^{\infty}(-1)^nq^ne^{[n]_qt}.
\tag13$$ Note that $\lim_{q\rightarrow 1}F_q(t)=\frac{2}{e^t +1}.$
By (11), we see that $\lim_{q\rightarrow 1}K_{n,q}=E_n.$ Hence,
$K_{n,q}$ seems to be the nice $q$-analogue of $E_n$.

\enddocument